\documentclass[a4paper,reqno]{amsart}
\pdfoutput=1

\usepackage{amsmath,amsfonts,amssymb,amsthm,stmaryrd,bbm,graphicx,mathtools,enumerate}
\usepackage{mathrsfs}
\usepackage[T1]{fontenc} 
\usepackage[latin1]{inputenc}
\usepackage[english]{babel}
\usepackage[tracking,spacing,kerning,babel]{microtype}
\usepackage{wasysym} 
\usepackage{color}
\usepackage{xcolor}
    	\usepackage{framed} 
\usepackage{wrapfig} 

\usepackage[final]{hyperref}   
\hypersetup{
    linktoc=page,
    linkcolor=red,          
    citecolor=blue,        
    filecolor=blue,      
    urlcolor=cyan,
   colorlinks=true           
}

\usepackage{lipsum} 

\usepackage{minitoc}
\usepackage{multicol}

\newtheoremstyle{mine}
{\baselineskip}
{\baselineskip}
{\itshape}
{
}
{\bfseries}
{.}
{.5em}
{#1 #2\ifx#3\relax\else~(#3)\fi}

\theoremstyle{mine}

\newtheorem{theorem}{Theorem}
\numberwithin{theorem}{section}
\newtheorem{corollary}[theorem]{Corollary}
\newtheorem{proposition}[theorem]{Proposition}
\newtheorem{lemma}[theorem]{Lemma}

\newtheorem{definition}[theorem]{Definition}
 
\numberwithin{equation}{section}

\theoremstyle{remark}
\newtheorem{remark}{Remark}

\newtheorem{question}{Question}

\colorlet{shadecolor}{blue!10}

\def\rm{\reversemarginpar}

\let\qed=\QED

\renewcommand{\epsilon}{\varepsilon}

\newcommand{\R}{\mathbb{R}}

\newcommand{\C}{\mathbb{C}}
\newcommand{\Z}{\mathbb{Z}}

\renewcommand{\S}{\mathbb{S}}

\def\Re{{\rm Re}\,}
\def\Tr{{\mathrm{Tr}}\,}

\def\P{\mathbb{P}} 
\def\E{\mathbb{E}} 
\def\md{\mid}

\def \eps {\epsilon}

\def\EFK#1#2#3{{\def\md{\bigm| } \E_{#1}^{\,#2}  \bigl[  #3 \bigr]}}

\def \p {{\partial}}

\def\<#1{\langle #1\rangle}

\definecolor{darkgreen}{rgb}{0,0.6,0.05}

\def\nn{\nonumber}
\def\bi{\begin{itemize}}  
\def\ei{\end{itemize}}
\def\bnum{\begin{enumerate}} 
\def\enum{\end{enumerate}}
\def\ni{\noindent}
\def\bf{\bfseries}

\newcommand{\Mod}[1]{\ (\mathrm{mod}\ #1)}

\def\Heis{\mathrm{Heis}}

\title[]{Continuous Symmetry Breaking along the Nishimori line}

\author{Christophe Garban}
\author{Thomas Spencer}

\address
{Université Claude Bernard Lyon 1, CNRS UMR 5208, Institut Camille Jordan, 69622 Villeurbanne, France \, and Institut Universitaire de France (IUF)}
\email{garban@math.univ-lyon1.fr}

\address{Institute for Advanced Study, 1 Einstein Drive, Princeton NJ 08540, USA}
\email{spencer@math.ias.edu}

\begin{document}

\maketitle

\begin{center}
{\em Dedicated to the memory of Freeman Dyson }
\end{center}

\begin{abstract}
We prove continuous symmetry breaking in three dimensions for a special class of disordered models  described by the Nishimori line. The spins take values in a group such as $\S^1$, $SU(n)$ or $SO(n)$.
Our proof is based on a theorem about group synchronization proved by Abbe, Massoulié, Montanari, Sly and Srivastava   \cite{abbe2018group}.  It also relies on a gauge transformation acting jointly on the disorder and the spin configurations due to Nishimori \cite{nishimori1981internal,georges1985exact}.  The proof does not use reflection positivity. The correlation inequalities of \cite{Messager1978correlation} imply symmetry breaking for the classical $XY$  model without disorder.
\end{abstract}


\section{Introduction}

\subsection{Context.} 

Nishimori introduced in \cite{nishimori1981internal} a one parameter family of random bond Ising models with a special quenched disorder on  edges, which is now called the {\em Nishimori disorder} or {\em Nishimori line}. 
 This special disorder was later discussed for more general symmetry classes by \cite{georges1985exact,ozeki1993phase, nishimori2001statistical, nishimori2002exact} and \cite{singer2011angular,wang2013exact}.  
The Nishimori line also has natural connections with Bayesian statistics and image processing, 
\cite{iba1999nishimori, nishimori2001statistical, tanaka2002statistical}. 

The Nishimori disorder is defined for all inverse temperatures $\beta$ and  it is associated with a special {\em gauge symmetry} which ensures that there is no singularity in the specific heat as $\beta$ varies. Furthermore, local energy correlations are  independent for all temperatures (see Lemma \ref{l.factor} below which will be a key step in our proof). Despite this lack of singularity, we shall see that spin correlations acquire a long range order at low temperature in 3 dimensions and above. 

Our proof of  long range order (symmetry breaking) relies heavily on \cite{abbe2018group} and Nishimori gauge invariance. It emphasizes the relation between phase transitions of special disordered $O(n)$ spin models and group synchronization in the presence of noise. This result was anticipated in \cite{abbe2018group}. 



The proof of continuous symmetry breaking in three dimensions is a challenging mathematical problem especially in case the group is not abelian, such as $SU(2)$. This article describes a new method to prove continuous symmetry breaking for spin systems with  Nishimori disorder. For spin models without disorder, the  proof of continuous symmetry breaking in three dimensions was established in \cite{frohlich1976infrared}. This proof relies on translation invariance and reflection positivity.  See also \cite{biskup2009reflection} for a review. 

Our results on long range order rely on the presence of the Nishimori disorder defined below. However, J\"urg Fr\"ohlich noted that for the $XY$ model, the long range order obtained on the Nishimori line implies long range order for the model without disorder. This result follows from a correlation inequality of Messager, Miracle-Sole and Pfister, \cite{Messager1978correlation}. The statement and proof of this inequality are given for completeness in the appendix. 

For the $O(2)$ symmetry class,  the works \cite{guth1980existence,frohlich1982massless,kennedy1986spontaneous} apply to more general cases without reflection positivity. Moreover, Ginibre inequalities can be used to prove long range order for disordered ferromagnetic models by comparison to the translation invariant case analyzed in \cite{frohlich1976infrared}.   In a series of papers,  Balaban (\cite{balaban1995low,balaban1996low,balaban1998low,balaban1998large}) developed  robust techniques to prove $O(n)$ symmetry breaking. It is possible that his analysis could be used to establish the results of this article. We are not aware of any other methods that would work in the non-abelian case for Nishimori disorder. Note that in one or two dimensions, symmetry breaking cannot occur by the Mermin-Wagner theorem \cite{mermin1966absence}.

One advantage of the present proof is that, in contrast to reflection positivity methods, it is mostly indifferent to the choice of domains in $\Z^d$ ($d\ge 3$)  as well as coupling constants  which may be spatially dependent. 




We now briefly describe the article by Abbe et al which is our main source of inspiration. They consider group elements $g_j$  with $j$ on the grid $\mathbb{Z}^3$.  Given  information  about $g_i^{-1}g_{j} +\,noise$ where $i$ and $j$ are adjacent pairs they prove that for small {\em noise}, information can be recovered about the long distance relative orientations $g_0^{-1}g_x$ for $x\in \Z^3$ far from 0. 
They provide an elegant and new reconstruction algorithm which is based on the notion of  {\em unpredictable paths} from \cite{benjamini1998unpredictable} and which will be of key use below.

In \cite{garban2020statistical}, the following related statistical reconstruction problem about the two dimensional Gaussian free field $\phi(j)$ has been considered: given $e^{i\alpha \phi(j)}$ (for some fixed  $\alpha \in \mathbb{R}$) when can we recover large scale information about $\phi$? I.e, instead of considering $g_i^{-1}g_j + noise$ as in \cite{abbe2018group}, the available data in \cite{garban2020statistical} is rather $g_i^{-1}g_j \Mod{2\pi/\alpha}$. This latter work is related to the Kosterlitz-Thouless transition \cite{kosterlitz1973ordering}.

This article is organized as follows. We first define the Nishimori disorder and state our main result for the disordered $XY$ model. Here the spins take values in the unit circle $\mathbb{S}^1$.  We  then formulate similar results for  disordered models with spins in $SU(2)$ and more general Lie groups. In section \ref{s.XY}, the proof of long range order is given for the $XY$ case. The proof begins with Lemma \ref{l.factor} which shows that the gauge invariance along Nishimori line implies independence of local energy correlations.   Theorem \ref{th.UP} from \cite{benjamini1998unpredictable} defines a measure on paths in $\mathbb{Z}^3$ with the property that two such paths rarely intersect. This theorem only holds in three or more dimensions. Exactly as in \cite{abbe2018group}, such paths are used to connect distant spins. Then  these two results are combined to prove long range order.  Section \ref{s.G} discusses the case when the spins take values in a Lie group. The proof follows as in the $XY$ case.  

A Nishimori disorder for the classical Heisenberg model when the spins take values in the sphere $\mathbb{S}^2$  is discussed in section \ref{s.Heis}. Note that in this case, for a given realization of the Nishimori disorder, the action is no longer invariant under $SO(3)$. In Section \ref{s.isoc}, we introduce a different Nishimori line in case where spins are in the unit sphere $\S^3 \equiv SU(2)$ for which we prove a symmetry breaking of right-isoclinic rotations. 
 Section \ref{s.final} comments on some generalizations and discusses open questions. This article concludes with an appendix proving the correlation inequality of Messager et al.


\subsection{Disordered $XY$ model and the Nishimori line.} 

We will start by introducing the relevant disorder in the special case of the disordered $XY$ model. We refer to \cite{nishimori1981internal,georges1985exact} for the more classical cases of random-bond Ising and Potts models and for a clear explanation why a specific line, the {\em Nishimori line} stands out among the possible disorders.  

\subsubsection{The Nishimori line.}
We fix a finite graph $\Lambda \subset \Z^d$.
The standard $XY$ model on $\Lambda$ has spins in $\S^1$ parametrized by angles $\theta(i)\, \in [0, 2\pi)$. At inverse temperature $\beta$  the Gibbs weight for free boundary conditions is proportional to 
\begin{align}\label{e.XY}
e^{\beta \sum_{i\sim j} \cos(\theta(i) - \theta(j))}\prod d\theta(i)\,.
\end{align}
The sum above ranges over all nearest neighbor pairs.  We will denote the Gibbs measure by $\<{\cdot}(\beta,\Lambda)$ and with those notations, the $XY$ spin-spin correlation will correspond to 
\begin{align}\label{e.SpinSpin}
\<{\cos (\theta(x) -\theta(y))}(\beta, \Lambda), \quad x,y \in \Lambda \subset\mathbb{Z}^d. 
\end{align}
If $\Lambda$ is symmetric enough and is equipped with periodic conditions, then when $d\geq 3$,  reflection positivity (RP) and infrared bounds (IRB) \cite{frohlich1976infrared} are used to prove that for $\beta$ large
there is long range order. (See also \cite{frohlich1982massless} which in this abelian setting avoids relying on reflection positivity). This means~\eqref{e.SpinSpin} is bounded below uniformly in $\Lambda$ and $x,y \in \Lambda$. 

We now introduce the following particular family of quenched disorders for this model. 
\begin{definition}[disordered $XY$ model]\label{d.disXY}
Let $\Lambda\subset \Z^d$ be a finite box and fix  $u>0$ and $\beta>0$. 
A {\em quenched disorder} $\omega$ will correspond to a family of  real random variables $\omega= (\omega_{i,j})_{(i,j)\in \vec E(\Lambda)}$ assigned to the oriented edges $\vec E(\Lambda)$ of $\Lambda$ and satisfying the following conditions:
\bnum
\item For any $e=(i,j)\in \vec E(\Lambda)$, $\omega_{i,j} = -\omega_{j,i}$. I.e. $\omega:=(\omega_{i,j})_{(i,j)\in \vec E(\Lambda)}$ may be viewed as a random $1$-form on $\Lambda$. 
\item For any two edges $e,e'$ which are not associated to the same undirected edge, $\omega_e$ is independent of $\omega_{e'}$. 
\item For each given edge $e=(i,j)$, the law of $\omega_{i,j}$ is supported on $[-\pi,\pi)$ and its density is given by the distribution
\begin{align}\label{e.rho}
\rho_u(\omega_{i,j}):= Z_{\rho_u}^{-1} e^{u \cos(\omega_{i,j})}
 \quad\text{ with } \quad  Z_{\rho_{u}}:= \int_{-\pi}^{\pi} e^{u \cos(\omega_{i,j}) } d \omega_{i,j}\,.
\end{align}
\enum
We will denote by $\P_{u}^{XY}$ and $\E_{u}^{XY}$ the probability measure and expectation corresponding to this quenched disorder $\omega$ on $[-\pi,\pi)^{\vec E(\Lambda)}$.  Notice that as $u\to \infty$, $\P_{u}^{XY}$ converges in law to the case of no disorder $\omega\equiv 0$. 

Given a fixed disorder $\omega$, we consider the following modified Gibbs weight
\footnote{Notice the factor $\tfrac 1 2$ in the modified Gibbs weight which was not present in the classical $XY$ model in~\eqref{e.XY}. This is due to that fact that  in this less symmetric case one needs to sum over oriented edges (equivalently one may also choose a prescribed direction for each edge and remove $\tfrac 1 2$). Both definitions match when $u\to \infty$.}
\begin{align}\label{e.disXY}
e^{\frac \beta 2 \sum_{(i,j)\in \vec E(\Lambda)} \cos(\theta(i) - \theta(j) + \omega_{i,j})}  \prod d\theta(i)\,.
\end{align}
The corresponding quenched partition function and expectation will be denoted by 
\begin{align}\label{}
Z_{\omega,\beta,\Lambda} \quad \text{   and   } \quad \<{\cdot}_{\omega, \beta, \Lambda} \,. 
\end{align}
\end{definition}
This gives us a two parameter  family $(u,\beta)$ of $XY$ models in random bond environment. The Nishimori line corresponds to the following special line which will satisfy extra integrability properties. 
\begin{definition}[Nishimori line]\label{d.nishXY}
The Nishimori line corresponds to the case where $u =\beta$. For any $\beta>0$, we will often work with its associated averaged quenched Gibbs measure
given by 
\begin{align*}\label{}
\EFK{\beta}{XY}{\<{\cdot}_{\omega,\beta,\Lambda}}\,.
\end{align*}
\end{definition}

\subsubsection{Adding a Random field.}
Given $h=\{h_j\}_{j\in \Lambda}\in \R^\Lambda$,
 we can also analyze the case when a random magnetic field $ \sum_j h_j\cos(\theta(j)+\psi(j))$ is present. Here $\psi(j)$
denote independent (not necessarily identically distributed) random phases. On the corresponding Nishimori line, the distribution of $\psi_j$ is  $Z_h^{-1}e^{\sum_j h_j\cos(\psi(j))}$ where $Z_h$ is the normalization factor. 
The  quenched expectation is now denoted by $\<{\cdot}_{\omega,\beta,\psi, h,\Lambda}$ and the expectation over the disorder is  $\mathbb {E}^{XY}_{\beta,h}$.

\subsubsection{Main result for the XY model on the Nishimori line.}
We now state our first main result in the case of a finite domain. (The corresponding statement in the non-abelian case will be given in Theorem \ref{th.disG}). 

\begin{theorem}\label{th.disXY}
	Let $d\geq 3$. There exist constants $\beta^*,a >0$  such that for  all $\beta\geq \beta^*$, all magnetic fields $h=\{h_j\}$, all finite $\Lambda \subset \Z^d$ and  all points $x,y \in \Lambda$ s.t. $B(\tfrac{x+y} 2, 2\|x-y\|_2) \subset \Lambda$, 
	\begin{align}\label{e.LXY}
	\EFK{\beta, h}{XY}{\<{\cos(\theta(x) - \theta(y))}_{\omega,\psi,\beta,h,\Lambda}} \geq 1 - \sqrt{\frac {a \log \beta} \beta}\,.\end{align} 
N.B. Continuous symmetry  holds when $h\equiv 0$ as the interaction is invariant under the global rotation  $\theta_j \rightarrow \theta_j + \alpha$.
\end{theorem}

\begin{remark}
	By the correlation inequality of \cite{Messager1978correlation} (see Appendix) when $h=0$, 
	$$\langle\cos(\theta(x) - \theta(y))\rangle_{\omega=0,\beta,\Lambda} \geq \langle\cos(\theta(x) - \theta(y))\rangle_{\omega,\beta,\Lambda}    $$ hence Theorem 1.3 implies long range order without disorder.
\end{remark}

\begin{remark}
By a standard high-temperature expansion, it is easy to check that for $h=0$ and $\beta$ small 
there exist $c,C>0$ such that, 
\begin{align}\label{}
|\<{\cos(\theta(x) - \theta(y))}_{\omega,\beta,\Lambda}| \leq C \exp(-c \|x-y\|_2)\,,
\end{align}
uniformly in the disorder $\omega$ and in the choice of $x,y\in \Lambda$. 
\end{remark}


\begin{remark}
We believe the lower bound~\eqref{e.LXY} 
is proportional to $\beta^{-1}$ as it is the case in the absence of disorder.
\end{remark}

\begin{remark}
Notice that  even at high $u\equiv \beta$ the disordered XY model is not ferromagnetic. This prevents us from using Ginibre's inequality to define the  infinite volume limit for the model. We will discuss infinite volume limits below in Section \ref{ss.IVL}. 
\end{remark}

\subsubsection{Dirichlet boundary conditions.}

The above definitions handle the case of {\em free boundary} conditions.
Let us define $0$-Dirichlet boundary conditions in our setting. We consider  any {\em boundary} $\p \Lambda \subset \Lambda$ (which may only consist in the classical interior boundary of $\Lambda$ but may also include, if desired, some bulk points). For any boundary site $i\in \p \Lambda$, fix $\theta(i):=0$ (east oriented). With this setup, the law on the quenched disorder remains unchanged while the (quenched) Gibbs weight becomes 
\begin{align}\label{e.disXY-D}
e^{\frac \beta 2 \sum_{(i,j)\in \vec E(\Lambda)} \cos(\theta(i) - \theta(j) + \omega_{i,j})}  \prod_{i\in \p \Lambda} \delta_0(\theta(i)) \prod_{i\in \Lambda\setminus \p \Lambda} d\theta(i)\,.
\end{align}
The corresponding quenched partition function, Gibbs measure and averaged quenched measure (on the Nishimori line) are denoted by 
\begin{align*}\label{}
Z^0_{\omega,\beta,\Lambda} \,, \quad \<{\cdot}^0_{\omega, \beta, \Lambda} \,, \quad \EFK{\beta}{XY}{ \<{\cdot}^0_{\omega, \beta, \Lambda} }\,. 
\end{align*}

Notice that such
Dirichlet boundary conditions are equivalent to applying a strong magnetic field $h_j \rightarrow \infty$ 
	for every $j\in \p \Lambda$. In this setting all the proof of long-range order goes through without a change (see Remark \ref{r.Dproof}). In particular this implies readily the following Corollary.
\begin{corollary}\label{th.disXYD}
With the same setup as in Theorem \ref{th.disXY}, the system acquires spontaneous magnetization for all $\beta \geq \beta^*$ in the following sense, let $\Lambda_n:=\{-n,\ldots,n\}^d \subset \Z^d$, then 
	\begin{align}\label{e.onepoint}
	\liminf_{n\to \infty} \EFK{\beta}{XY}{\<{\cos(\theta(0))}^0_{\omega,\beta,\Lambda_n}} \geq 1-
	\sqrt{\frac {a \log \beta} \beta}\,.
	\end{align}
\end{corollary}


\begin{remark}\label{}
It would be tempting at this point to recover more general boundary conditions such as {\em Dobrushin boundary conditions} by letting $h_j\to -\infty$ on one side of the boundary. We wish to emphasize here that this trick would not produce true Dobrushin boundadry conditions. Indeed, the random phase $\psi(j)$ on points $j\in \p\Lambda$ with $h_j \to -\infty$ would then need to be sampled according to $e^{h_j \cos(\psi(j))}$ which converges to $\delta_{\pi}$ and therefore produces effectively the same effect as setting $h_j \to +\infty$. 
%
\end{remark}

\subsection{Disordered $SU(2)$ spin model (and other such Lie groups).}\label{ss.G}
As mentioned earlier, the main interest of this work is that it provides a robust way, using statistical reconstruction techniques, to deal with spin systems with non-abelian symmetry. We introduce the relevant such spin models in this Section. Following \cite{abbe2018group} which deals with {\em group synchronization}, the  spin systems we will be able to analyze will carry at any vertex $x$,  a group elements $U_x$ in some given compact matrix Lie group $G$, for example $G=$
\begin{align}\label{e.list}
SU(2) \quad\, SO(3) \quad SU(n) \quad U(n) \quad \text{ etc.} 
\end{align}

Before introducing a quenched disorder, let us first introduce our spin system on $G^{\Lambda}$ where $G$ is one of the above groups and $\Lambda$ is a finite graph in $\Z^d$. The spin configurations will be denoted $U=(U(x))_{x\in \Lambda}$ where each $U(x)\in G$. At inverse temperature $\beta$  the Gibbs weight for free boundary conditions is proportional to 
\begin{align}\label{e.XYG}
e^{\beta \sum_{i\sim j} \Re (\Tr U^*(i) U(j) ) }\prod d\mu_G(U(i))\,.
\end{align}
The sum above ranges over all (non-oriented) nearest neighbor pairs and $\mu_G$ denotes the normalized Haar measure on $G=SU(2), \, SO(3)\, $ etc. 

%

By analogy with the case of disordered $XY$ model, we define the following disordered model whose quenched disorder will still be denoted as $\omega$.

\begin{definition}[disordered non-abelian models]\label{d.disG}
Let $\Lambda\subset \Z^d$ be a finite box and fix  $u>0$ and $\beta>0$. 
Fix the group $G$ to be one of the above groups ($SU(2),\, SO(3)\,$ etc.). 
A {\em quenched disorder} $\omega$ will now correspond to a family of $G$-valued matrices $\omega= (\Omega_{i,j})_{(i,j)\in \vec E(\Lambda)}$ assigned to the oriented edges $\vec E(\Lambda)$ of $\Lambda$ and satisfying to the following conditions
\bnum
\item For any $e=(i,j)\in \vec E(\Lambda)$, $\Omega_{i,j} = \Omega^*_{j,i}$. 
\item For any two edges $e,e'$ which are not associated to the same undirected edge, $\Omega_e$ is independent of $\Omega_{e'}$. 
\item For each given edge $e=(i,j)$, the law of $\Omega_{i,j}$ has the following distribution
\begin{align}\label{e.rhoG}
d \rho_u(\Omega_{i,j}):= Z_{\rho_u}^{-1} e^{u \Re(\Tr \Omega_{i,j})} d \mu_G(\Omega_{i,j})
 \quad\text{ with } \quad  Z_{\rho_{u}}:= \int_G e^{u \Re( \Tr \Omega_{i,j} ) } d \mu_G(\Omega_{i,j})\,.
\end{align}
\enum
We will denote by $\P_{u}^{G}$ and $\E_{u}^{G}$ the probability measure and expectation corresponding to this quenched disorder $\omega$ on $G^{\vec E(\Lambda)}$.  Notice that as $u\to \infty$, $\P_{u}^G$ converges in law to the deterministic disorder $\omega\equiv \mathrm{Id}_G$. 

Given a fixed disorder $\omega$, we consider the following modified Gibbs weight
\begin{align}\label{e.disG}
e^{\frac \beta 2 \sum_{(i,j)\in \vec E(\Lambda)} \Re(\Tr U^*(i) \Omega_{i,j} U(j) )}  \prod d\mu_G(U(i))\,.
\end{align}
The corresponding quenched partition function and expectation will be denoted (with a slight abuse of notations as we do not distinguish these notations from the $XY$ case) by 
\begin{align}\label{}
Z_{\omega,\beta,\Lambda} \quad \text{   and   } \quad \<{\cdot}_{\omega, \beta, \Lambda} \,. 
\end{align}

As for the disordered XY model, the \textbf{Nishimori line} of this model will correspond to the line $u\equiv \beta$, i.e. to the averaged quenched measures
\begin{align*}\label{}
\EFK{\beta}{G}{\<{\cdot}_{\omega,\beta,\Lambda}}\,.
\end{align*}
\end{definition}

\begin{remark}\label{r.GlobInv}
For a fixed disorder $\omega=(\Omega_{i,j})_{(i,j)\in \vec E(\Lambda)}$, note that the action in~\eqref{e.XYG} is invariant under a global $G$-rotation because of the trace.
\end{remark}

\ni
\textbf{Boundary conditions.} As in the case of the disordered $XY$ model, we will be able to treat the following two boundary conditions along the Nishimori line. 
\bi
\item {\em Free boundary conditions}. This corresponds to the above Definition. 
\item {\em Dirichlet boundary conditions.} For any subset $\p \Lambda \subset \Lambda$, we set the spins $\{U(i)\}_{i\in \p \Lambda}$ to be equal to $\mathrm{Id}_G$ and the law of the disorder on all oriented edges in $\vec E(\Lambda)$ is exactly the same as for free boundary conditions.  
\ei
We may now state our main result in the non-abelian setting where, as before, we fix $G$ to be any given compact matrix Lie group $G$ (for example $G=
SU(2),\, SO(3),\, SU(n),\,  U(n)$  etc).

\begin{theorem}\label{th.disG}
Let $d\geq 3$. There exist constants $\beta^*,a >0$  such that for any  $\beta\geq \beta^*$, any  $\Lambda \subset \Z^d$ and any points $x,y \in \Lambda$ s.t. $B(\tfrac{x+y} 2, 2\|x-y\|_2) \subset \Lambda$, 
\begin{align}\label{}
\EFK{\beta}{G}{\<{\Re( \Tr U^*(x) U(y))}_{\omega,\beta,\Lambda}} \geq \big(1 - \sqrt{\frac {a \log \beta} \beta}\big) \Tr \mathrm{Id}_G \,.
\end{align}
Furthermore, the system acquires spontaneous magnetization for all $\beta \geq \beta^*$ in the following sense, let $\Lambda_n:=\{-n,\ldots,n\}^d \subset \Z^d$, then 
\begin{align}\label{}
\liminf_{n\to \infty} \EFK{\beta}{G}{\<{\Re( \Tr U(0))}^0_{\omega,\beta,\Lambda_n}} \geq 
\big(1 - \sqrt{\frac {a \log \beta} \beta}\big) \Tr \mathrm{Id}_G \,.
\end{align} 
\end{theorem}

\begin{remark}\label{}
In the case of the classical spin $O(n)$ model (including the classical Heisenberg model) where spins take their values in the unit sphere $\S^{n-1}$, the fact we cannot rely on an underlying group structure for the displacements $U_i^{-1} U_j$ raises some difficulty. 
In Section \ref{s.Heis} we will introduce a model of (anisotropic) quenched disorder for the classical Heisenberg model for which we will be able to prove a phase transition in $d\geq 3$.
(See Theorem \ref{th.heis}). In the special case of the spin $O(4)$-model we introduce a more isotropic Nishimori line for which we can prove a symmetry breaking of right-isoclinic rotations. This is discussed in Section \ref{s.isoc}. 
In order to lighten the introduction, we postpone this analysis of classical spin $O(n)$ models  to Sections \ref{s.Heis} and \ref{s.isoc}. 
\end{remark}



\subsection{Infinite volume limits and spontaneous magnetization.}\label{ss.IVL} 
 By compactness of $\S^1$ and $G=SU(2), \, SO(3)$ etc., the existence of infinite volume limits under the averaged quenched measure $\EFK{\beta}{G}{\<{\cdot}_{\omega,\beta,\Lambda}}$ is straightforward. Yet the question of uniqueness in this averaged case is much less clear. Even the existence of quenched infinite volume limits which are measurable w.r.t the disorder $\omega$ turns out to be rather subtle as we shall see below.  
  In fact, already in the classical case (without disorder), the uniqueness of Gibbs measures on infinite lattices is not known for non-abelian groups $G$. This is due to the lack of Ginibre's inequality in these cases (see Ginibre \cite{ginibre1970general}).
Finally, let us stress that if one samples a Nishimori disorder $\omega\sim \E_\beta^{XY}$ on the whole lattice $\Z^d$, then our proof of Theorem \ref{th.disXY} does not imply for example that
$\EFK{\beta}{XY}{\liminf_{n\to \infty} \<{\cos(\theta_0)}^0_{\omega,\beta,\Lambda_n}} \geq 1 - \sqrt{\frac {a \log \beta} \beta}$ holds.

All these facts show that some care is needed when dealing with infinite volume limits. 

For the construction of measurable quenched infinite volume limits, one proceeds as follows: consider the sequence of couplings $(\omega, \<{\cdot}^0_{\omega,\beta,\Lambda_n})_n$, where $\omega$ may be viewed either as a random environment on $\Lambda_n$ or directly on the full lattice $\Z^d$. Then one can argue by compactness that there exist subsequential limits in law $(\omega, \widetilde{\<{\cdot}}^0_{\md \omega, \beta, \Z^d})$. The subtlety here is that the Gibbs measure $\widetilde{\<{\cdot}}^{0}_{\md \omega, \beta,\Z^d}$ may no longer be a deterministic function of $\omega$ and it may only  be a random Gibbs measure conditionally on $\omega$ (this is why we denote it with a subscript $|\omega$). We use the notation $\widetilde{\<{\cdot}}$ to stress here that we are taking a subsequential scaling limit. In this slightly weaker quenched sense, Theorem \ref{th.disG} readily implies the following spontaneous magnetization result in infinite volume.

\begin{corollary}[Spontaneous magnetization]\label{c.SM}
Fix $\beta>0$ and $G$ any Lie group from the list~\eqref{e.list}. 
Let $\omega\sim \E_{\beta}^G$ be a Nishimori disorder on the whole lattice $\Z^d$ and let $(\omega,\widetilde{\<{\cdot}}_{\md\omega, \beta, \Z^d})$ be any subsequential infinite volume limit as outlined above. Then, with the same notations as in Theorem \ref{th.disG} we have for any $\beta\geq\beta^*$,
\begin{align*}\label{}
\EFK{\beta}{G}{\widetilde{\<{\Re( \Tr U(0))}}^0_{|\omega,\beta,\Z^d}} \geq 
\big(1 - \sqrt{\frac {a \log \beta} \beta}\big) \Tr \mathrm{Id}_G\,.
\end{align*}
\end{corollary}

\subsection{Acknowledgements.}
We wish to thank Roland Bauerschmidt, Christophe Sabot and Avelio Sepúlveda for useful discussions. We thank J\"urg Fr\"ohlich for explaining how our results on the Nishimori line imply long range order for the $XY$ model without disorder. Finally we thank the anonymous referee for a very careful reading of the manuscript. The research of C.G. is supported by the ERC grant LiKo 676999 and the Institut Universitaire de France.


\section{Symmetry breaking for disordered XY model}\label{s.XY}

We start with the following lemma which reveals the significance of the Nishimori line. Its proof will use crucially a certain {\em gauge transformation} which does not seem to play an important role in \cite{abbe2018group} but played a central role in the original paper \cite{nishimori1981internal} by Nishimori on Ising and Potts models. See also \cite{georges1985exact} where such gauge transformations have been used to compute explicitly several quantities (such as the averaged quenched internal energy) on the Nishimori line. 

Let $\{ f_{i,j}\}_{(ij) \in E(\Lambda)}$ be smooth and periodic test functions and let $\prod_{(ij) }$ denote the product of the edges in $\Lambda$. The integrals below are product integrals which  range over $[0, 2\pi]$.
 
\begin{lemma}\label{l.factor}
For any finite domain $\Lambda \subset \Z^d$, any inverse temperature $\beta$ and any magnetic fields $\{h_j\}$,  then under the Nishimori-line averaged quenched measure $$\EFK{\beta,h}{XY}{\<{\prod_{ij} f_{i,j}(\theta(i) - \theta(j) +\omega_{i,j})}_{\omega,\psi,\beta, h,\Lambda}} = \mathbb{E}_{\beta}[\prod_{ij} f_{i,j}(\omega)] $$
where $\mathbb{E}_{\beta}$ is the disorder average over $\omega$. The right side is independent of $h$.
\end{lemma}

\noindent{\it Proof.} 
To avoid dealing with the two possible orientations of each edge, let us choose an arbitrary direction for each unoriented edge $e$. We shall denote by $E(\Lambda)\subset \vec E(\Lambda)$ this subset.

By definition
\begin{align}\label{e.factor1}
&\EFK{\beta, h}{XY}{\<{\prod_{(ij) }  f_{i,j} (\theta(i) - \theta(j) + \omega_{i,j})}_{\beta, h,\Lambda}} \\
&  =\mathbb{E}_{\beta, h}  {\frac 1 {Z_{\omega,\psi,\beta,h,\Lambda}}  \int  \prod_{(ij)}  f_{i,j}(\theta(i)-\theta(j)+\omega_{i,j})  e^{\beta \cos(\theta(i) - \theta(j) + \omega_{i,j})}
	 \prod_j e^{ h_j\cos(\theta(j)+\psi(j))}  d\theta(j)}\,. \nn
\end{align}
Here $\mathbb{E}_{\beta, h}$ is the disorder expectation in $\omega, \psi$ given by 
\[
Z^{-1}_{\beta } \prod_{(ij)}e^{\beta\cos(\omega_{i,j})}\,Z^{-1}_h \prod_j e^{h_j\cos(\psi(j))}\,.
\]
Now, for any fixed deterministic field $\{\phi(i)\}_{i\in \Lambda}$ in $[0,2\pi)$, we make the following change of variables
\begin{align}\label{e.GT}
\begin{cases}
& \theta(i) \rightarrow \theta(i) - \phi(i)  \\
& \omega_{i,j} \rightarrow \omega_{i,j} +\phi(i) - \phi(j),\quad \psi(j) \rightarrow \psi(j) +\phi(j)  \,. 
\end{cases}
\end{align}
Note that this change of variables is chosen so that it does not affect the numerator on the right hand side of (2.1) defined by
\begin{align*}\label{}
\mathbf A:= \int  \prod_{(ij)}  f_{i,j}(\theta(i)-\theta(j)+\omega_{i,j})  e^{\beta \cos(\theta(i) - \theta(j) + \omega_{i,j})}\prod_j e^{ h_j\cos(\theta(j)+\psi(j))}  d\theta(j)\,.
\end{align*}
The only effect of this change of variables is to shift the weight of the disorder in $\mathbb{E}_{\beta,h}$. Namely, after the change of variables~\eqref{e.factor1}, we obtain  that for any prescribed field $\{\phi(j)\}_{i\in \Lambda}$, 
\begin{align}\label{e.factor2}
&\EFK{\beta, h}{XY}{\<{\prod_{(ij)}  f_{i,j} (\theta(i) - \theta(j) + \omega_{i,j})}_{\omega,\beta,h,\Lambda}} \nn \\
& =   \frac 1 {Z_{\beta}Z_h} \int  \frac {\mathbf A } {Z_{\omega,\psi,\beta,h,\Lambda}}   \prod_{(ij)} e^{\beta \cos(\omega_{i,j}+\phi(i) - \phi(j) )} d\omega_{i,j}\prod_j e^{ h_j\cos(\phi(j)+\psi(j))}  d\psi(j) \,.
\end{align} 

Now the key observation behind the  factorization on the Nishimori line is that since  (2.3) does not depend on the choice of field $\{\phi(i)\}_{i\in \Lambda}$, one may as well integrate the expression (2.3) over uniformly chosen $\phi(i)\in[0,2\pi)$. After integrating over $\phi$ note that the numerator  $$\int \prod_{(ij)} e^{\beta \cos(\omega_{i,j}+\phi(i) - \phi(j) )} \prod_j e^{h_j\cos(\phi_j +\omega_j)} d\phi(j)$$ cancels  $Z_{\omega,\psi,\beta,h,\Lambda}$. This cancelation enables  us to explicitly calculate the integral: 
\begin{align}\label{e.factor2}
&\EFK{\beta,h}{XY}{\<{\prod_{(ij) }  f_{i,j} (\theta(i) - \theta(j) + \omega_{i,j})}_{\omega, \psi,\beta,h,\Lambda}} \nn \\
& =   \frac 1 {Z_{\beta}\,Z_{h}} \int   
\mathbf A  \prod_{(ij)}  d\omega_{i,j} \, \prod_j d\psi(j)  \nn \\
& =  \frac 1 {Z_{\beta}\,Z_{h}}\int     
  \prod_{(ij)}  f_{i,j}(\theta(i)-\theta(j)+\omega_{i,j})   e^{\beta \cos(\theta(i) - \theta(j) + \omega_{i,j})} d\omega_{i,j}\prod_j e^{h_j\cos(\theta(j)+\psi(j)})d\psi(j) d\theta(j)\nn \\
 &  =  \frac 1{Z_{\beta}}     
 \int 
  \prod_{(ij)}  f_{i,j}(\theta(i)-\theta(j)+\omega_{i,j})  e^{\beta \cos(\theta(i) - \theta(j) + \omega_{i,j})}d\omega_{i,j}\prod_j d\theta(j) \nn \\
 & = \mathbb{E}_{\beta} \prod_{ij}f_{i,j}(\omega_{ij}) \,.
\end{align} 
 Note that we have normalized so that the integral over $\prod_j d\theta_j$ is one.  \qed

\medskip

We will derive from the above Lemma  a surprising property on the quenched internal energy along the Nishimori-line which is well-known and was established earlier on (see \cite{nishimori1981internal,georges1985exact,nishimori2001statistical,nishimori2002exact,jacobsen2002phase}). 

Given a domain $\Lambda$ and a disorder $\omega$, the quenched internal energy corresponds to 
\begin{align*}\label{}
\<{E}_{\omega,\beta,\Lambda} & := -\frac {\p \log Z_{\omega,\beta,\Lambda}} {\p \beta} \\
& = -\frac
{\int   \frac 1 2 \sum_{(i,j)\in \vec E(\Lambda)} \cos(\theta(i) - \theta(j) + \omega_{i,j})  e^{\frac \beta 2 \sum_{(i,j)\in \vec E(\Lambda)} \cos(\theta(i) - \theta(j) + \omega_{i,j})}   d\theta}
{\int e^{\frac \beta 2 \sum_{(i,j)\in \vec E(\Lambda)} \cos(\theta(i) - \theta(j) + \omega_{i,j})}  d\theta} 
\end{align*}



\begin{corollary}[\cite{nishimori2001statistical,nishimori2002exact}]\label{c.FE}
The (averaged) quenched internal energy on the Nishimori line can be explicitly computed. For any finite $\Lambda\subset \Z^d$, it is given by 
\begin{align*}\label{}
\EFK{\beta}{XY}{\<{E}_{\omega,\beta,\Lambda}} = - \frac {\int_\pi^\pi \cos(w) e^{\beta \cos w} dw} {Z_{\rho_\beta}}  |E(\Lambda)|\,,
\end{align*}
where $|E(\Lambda)|$ denotes the number of (non-oriented) edges in $\Lambda$. 

\end{corollary}

\ni
{\em Proof.}
It is an immediate corrollary of Lemma \ref{l.factor}.
\qed

\begin{remark}\label{}
Lemma \ref{l.factor} also implies that local energies at different edges factor
 under the averaged quenched measure (another related specificity of the Nishimori-line). This implies in turn that there is no singularity for the specific heat as $\beta$ varies. 
\end{remark}

As we shall explain further below, the identity below indicates that the Nishimori line does not enter the spin glass phase.  This identity was first proved in this Ising case in \cite{nishimori1981internal} and there are non-Abelian versions proved in \cite{ozeki1993phase,georges1985exact}.

\begin{corollary}
For $h=0$ and all $\beta$ the quenched expectation of the spin correlation is positive
	$$\mathbb{E}_{\beta}^{XY}[\<{e^{i(\theta_0-\theta_x)}}_{\omega,\beta,\Lambda}] = \mathbb{E}_{\beta}^{XY}[|\<{e^{i(\theta_0-\theta_x)}}_{\omega,\beta,\Lambda}|^2]  \ge 0. 
	$$ 
\end{corollary}
\ni
{\em Proof.}
Under the change of variables~\eqref{e.GT} and integrating over the $(\phi_i)$ variables, the obervable on the left-hand-side becomes $$\int\prod d\omega_{j,j'}Z_{\rho_{\beta}}^{-1}Z(\omega)\, \langle e^{i(\theta_0-\theta_x)}\rangle_{\omega,\beta,\Lambda}\langle e^{-i(\phi_0-\phi_x)}\rangle_{\omega,\beta,\Lambda}$$ where $Z_{\rho_{\beta}}$ is the normalization of the disorder and 
\begin{align*}
Z(\omega)= \int e^{\beta \sum_{jj'}\cos({ \hat \phi_j}-\hat \phi_{j'}+ \omega_{j,j'} )}\prod d{\hat \phi_j}\,.
\end{align*}
We now make another change of variables $\omega_{j,j'} \rightarrow \hat \omega_{j,j'} = \omega_{j,j'} + \hat{\phi_j} - \hat{\phi_{j'}} $. Notice that under this change of variables, one has 
\[
\langle e^{i(\theta_0-\theta_x)}\rangle_{\hat \omega,\beta,\Lambda}\langle e^{-i(\phi_0-\phi_x)}\rangle_{\hat \omega,\beta,\Lambda}
= 
\langle e^{i(\theta_0-\theta_x)}\rangle_{\omega,\beta,\Lambda}\langle e^{-i(\phi_0-\phi_x)}\rangle_{\omega,\beta,\Lambda}\,.
\]
Indeed the invariance of this product of expactations under $\omega \to \hat \omega$ follows from the fact that the observables are complex conjugates of each other. By Fubini and first integrating over $d\omega_{j,j'}$ we obtain 
\begin{align*}\label{}
& \mathbb{E}_{\beta}^{XY}[\<{e^{i(\theta_0-\theta_x)}}_{\omega,\beta,\Lambda}] \\ 
& = 
\int \prod d\hat \phi_j \int \prod d\omega_{j,j'} e^{\beta \cos(\omega_{j,j'})}Z_{\rho_{\beta}}^{-1}\, \langle e^{i(\theta_0-\theta_x)}\rangle_{\omega,\beta,\Lambda}\langle e^{-i(\phi_0-\phi_x)}\rangle_{\omega,\beta,\Lambda} \\
& = \int \prod d\hat \phi_j  \,\,  \mathbb{E}_{\beta}^{XY}[|\<{e^{i(\theta_0-\theta_x)}}_{\omega,\beta,\Lambda}|^2]  =  \mathbb{E}_{\beta}^{XY}[|\<{e^{i(\theta_0-\theta_x)}}_{\omega,\beta,\Lambda}|^2] \,. 
\end{align*}
 Moreover, using  inversions $\theta \rightarrow - \theta , \, \omega \rightarrow - \omega $ the expectations above equal  $\mathbb{E}_{\beta}^{XY}[\cos(\theta_0-\theta_x)] $.  \qed

\medskip

The reason why this corollary is a strong indication that the Nishimori line does not enter the {\em spin glass} phase is as follows: the left side of the identity is expected to go to zero in the SG phase as $x\to \infty$ while the right side is expected to remain bounded away from zero as $x\to \infty$. 
Moreover the Nishimori line is ``expected'' to pass through a multicritical point at the boundary of paramagnetic, ordered and spin glass phase. See \cite{nishimori1981internal,le1988location,ozeki1993phase, nishimori2002exact}.

\smallskip
The next Theorem
on the so-called {\em unpredictable paths} due to Benjamini, Pemantle and Peres (\cite{benjamini1998unpredictable}) plays a key role in \cite{abbe2018group} and will also be a key ingredient of our proof.  Their result concerns a probability measure  $\mu$ supported on infinite increasing paths in $\Z^3$ formed by sums of $(1,0,0), (0,1,0), (0,0,1)$ and starting at $0$. 


\begin{theorem}[\cite{benjamini1998unpredictable})]\label{th.UP} In three dimensions there exists a probability measure $\mu$ supported on increasing paths which satisfies the following 
{\em exponential intersection tails (EIT)} property. There exist constants $C,\alpha>0$ such that for any $k\geq 1$, 
\begin{align}\label{e.10}
(\mu \times \mu) \Big[(p_1,p_2)\,\, \text {such that}, \,\,|p_1 \cap p_2| \ge k \Big] \le C e^{-\alpha k}\,.
\end{align}
Here $|p_1 \cap p_2|$ denotes the number of common edges.
\end{theorem} 
	
The analog of this statement in $\Z^d, d\geq 4$ is much easier as the uniform measure on increasing paths satisfies this EIT property when $d\geq4$. In $d=3$ uniformly chosen paths do not satisfy the EIT  property. This is why unpredictable paths have been invented in \cite{benjamini1998unpredictable}.   
See also the papers by H\"aggstr\"om and Mossel and Hoffman \cite{haggstrom1998nearest, hoffman1998unpredictable}.

\medskip Inspired by Abbe et al \cite{abbe2018group} we now combine Lemma \ref{l.factor} and Theorem  \ref{th.UP} to prove our main theorem on the disordered $XY$ model. 

\smallskip
\noindent{\it Proof of Theorem \ref{th.disXY}.}
For simplicity, we shall stick to the case where $d=3$ and $h\equiv 0$, the case of higher dimensions being  easier as one does not need the construction of  unpredictable paths from Theorem \ref{th.UP}.
\medskip

\ni
\textbf{Step 1. Two point function between 0 and $\vec n =(n,n,n)$.}

Let $n\geq 1$. We assume here that the box $D_n=\{0,\ldots,n\}^3 \subset \Lambda$. Note that $D_n$ contains both $0$ and $\vec n$.

Define
\begin{align}\label{e.lambda}
\lambda = \lambda(\beta):= \EFK{\rho_\beta}{}{e^{i \omega}} = \frac {\int_{-\pi}^\pi \cos(\omega) e^{\beta \cos(\omega)} d\omega}{ \int_{-\pi}^\pi e^{\beta \cos(\omega)} d\omega}\,.
\end{align}
On can check that for large $\beta$, 
\begin{align}\label{e.lambda2}
\lambda = \lambda(\beta)= 1-  \frac 1 {2\beta} +o(\beta^{-1})\,.
\end{align}
For any $(i,j)\in \vec E(\Lambda)$, let 
\begin{align}\label{e.12}
Y_{i,j}:= e^{i(\theta(i)-\theta(j) +\omega_{i,j})}\,.
\end{align}
Lemma \ref{l.factor} readily implies that if $p$ is any non-intersecting path going from $0$ to $\vec n$, then 
\begin{align*}\label{}
\EFK{\beta}{XY}{\<{\prod_{(ij) \in p} Y_{i,j}}_{\omega,\beta,\Lambda}}
& = \EFK{\beta}{XY}{\<{e^{i(\theta(0) - \theta(\vec n))} \prod_{(ij) \in p} e^{i \omega_{i,j}}}_{\omega,\beta,\Lambda}} \\
&\hskip -3 cm  = \prod_{(ij)\in p} \EFK{\beta}{XY}{\<{Y_{i,j}}_{\omega,\beta,\Lambda}} \,\,\, \text{ (using the independance from Lemma \ref{l.factor})} \\
&\hskip -3 cm   = \lambda^{|p|}\,,
\end{align*}
where $|p|$ is the length of the path $p$. We obtain from these equalities the following identity which holds for any given simple path $p$ from $0$ to $\vec n$. 
\begin{align}\label{e.KI}
1 & = \EFK{\beta}{XY}{\<{e^{i(\theta(0) - \theta(\vec n))} \left(  \frac 1 {\lambda^{|p|}} \prod_{(ij) \in p} e^{i \omega_{i,j}} \right)}_{\omega,\beta,\Lambda}} \nn \\
& = \EFK{\beta}{XY}{\<{e^{i(\theta(0) - \theta(\vec n))}}_{\omega,\beta,\Lambda} \left(  \frac 1 {\lambda^{|p|}} \prod_{(ij) \in p} e^{i \omega_{i,j}} \right)}
\end{align}

As in \cite{abbe2018group}, we will average this identity over a suitably chosen probability measure on random paths $p$ from $0$ to $\vec n$. 
Following \cite{abbe2018group}, if $\mu$ is the measure on increasing paths from Theorem \ref{th.UP}, then for any $n\geq1$, by reflecting the paths $\gamma\sim \mu$ under the hyperplane $\{(x_1,x_2,x_3)\in \R^3,  x_1+x_2+x_3 = 3n/2\}$, we can define probability measures $\mu_n$ on simple paths $p$ from $0$ to $\vec n$ which satisfy the EIT property uniformly in $n$. We will still call these paths $p$ increasing paths. 
 
Let us define the following random variable which is measurable w.r.t. the Nishimori disorder $\omega$:
\begin{align}\label{e.R}
\mathbf R(\omega) = \mathbf R_{0,\vec n}(\omega):=  \frac 1 {\lambda^{3n}} \EFK{\mu_n}{}{\prod_{(ij) \in p} e^{i \omega_{i,j}}}\,. 
\end{align}
(We used here that all paths $p\sim \mu_n$ are such that $|p|=3n$). 
By averaging the identity~\eqref{e.KI} w.r.t to the measure on paths $\mu_n$, we obtain 
\begin{align}\label{e.great}
\EFK{\beta}{XY}{\<{e^{i(\theta(0) - \theta(\vec n))}}_{\omega,\beta,\Lambda} \; \mathbf R(\omega)} & =1 \,. 
\end{align}

Notice that by definition, $\mathbf R(\omega)$ is such that $\EFK{\beta}{XY}{\mathbf R(\omega)} =1$. Theorem \ref{th.disXY} will follow from the following second moment estimate on $\mathbf R(\omega)$ (whose proof is delayed to the end of this proof). 
\begin{lemma}[\cite{abbe2018group}]\label{l.R}
There exist $\beta^*,a>0$ such that for all $\beta \geq \beta^*$. 
\begin{align*}\label{}
\EFK{\beta}{XY}{|\mathbf R(\omega)|^2} \leq 1 + \frac {a \log \beta} \beta\,.
\end{align*}
\end{lemma}

A lower-bound on the two-point function between $0$ and $\vec n=(n,n,n)$ easily follows from this estimate. Indeed starting with the identity~\eqref{e.great}, we have
\begin{align*}\label{}
1& = \EFK{\beta}{XY}{\<{e^{i(\theta(0) - \theta(\vec n))}}_{\omega,\beta,\Lambda} \; \mathbf R(\omega)} \\
& = \EFK{\beta}{XY}{\<{e^{i(\theta(0) - \theta(\vec n))}}_{\omega,\beta,\Lambda}} 
+ \EFK{\beta}{XY}{\<{e^{i(\theta(0) - \theta(\vec n))}}_{\omega,\beta,\Lambda} \left( \mathbf R(\omega) - 1 \right)}\,.
\end{align*}
Thus,
\begin{align*}\label{}
1 - \EFK{\beta}{XY}{\<{\cos(\theta(0) - \theta(\vec n))}} & \leq 
|1- \EFK{\beta}{XY}{\<{e^{i(\theta(0) - \theta(\vec n))}}_{\omega,\beta,\Lambda}}| \\
& \leq   
\EFK{\beta}{XY}{\left|\mathbf R(\omega) - 1 \right|} \\
& \leq  \EFK{\beta}{XY}{|\mathbf R(\omega)|^2 -1}^{1/2} \,\, \text{  since } \EFK{\beta}{XY}{\mathbf R(\omega)} =1\\
&  \leq \sqrt{\frac {a \log \beta} \beta}\,,
\end{align*}
by Lemma \ref{l.R}. This concludes the proof when $x=0$ and $y=\vec n$. 
%

\begin{figure}[!htp]
\begin{center}
\includegraphics[width=\textwidth]{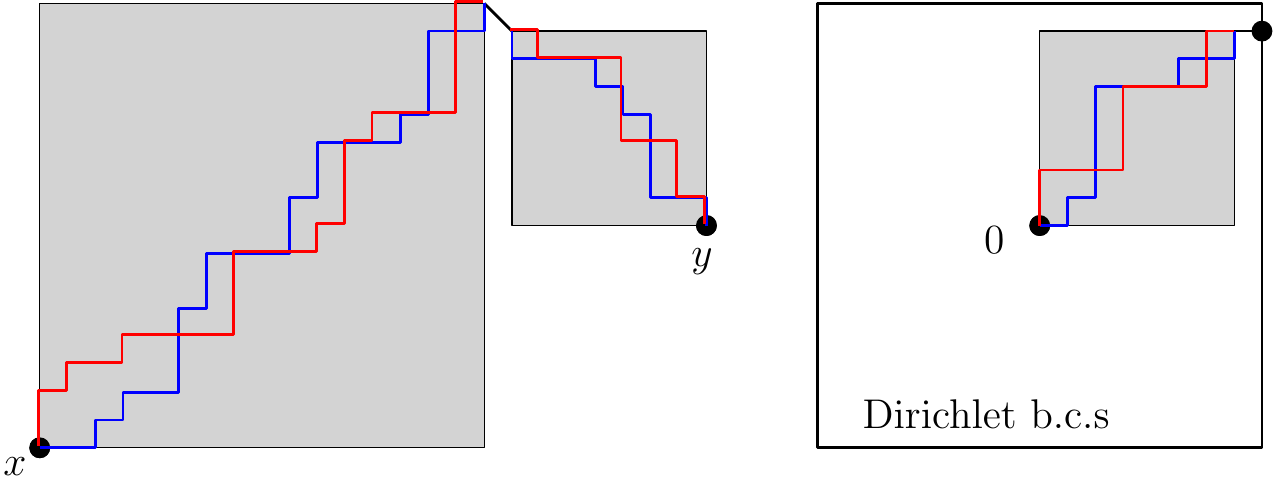}
\end{center}
\caption{Family of directed paths used from $x$ to $y$ and in the presence of Dirichlet boundary conditions.}\label{f.sailing}
\end{figure}

\medskip
\medskip
\ni
\textbf{Step 2. Two point function between $x$ and $y$ when $B(\tfrac{x+y} 2, 2\|x-y\|_2) \subset \Lambda$.}

Similarly as in \cite{abbe2018group}, the point is to notice that in $\Z^3$, one may travel from any $x$ to $y$ using at most 4 directed cones as in case 1 and which are non-intersecting. See Figure \ref{f.sailing} for an illustration in $\Z^2$. 
By translating the points and the domain, we may assume $x=0$ and $y=(a,b,c)$. WLOG, we may also assume $a\leq b \leq c$. Along the direction $(1,1,1)$, one may first go to $(a,a,a)$. Then, keeping the direction $(1,1,1)$ for a time $\tfrac{b-a} 2$ and then $(-1,1,1)$ for a time $\tfrac{b-a} 2$, we end up in $(a,b,b)$. In the remaining time $c-b$, we follow half-way $(1,1,1)$ and then $(-1,-1,1)$ to end up in $(a,b,c)$. To make this description more precise, $a,b,c$ should be assumed to be even here, otherwise we go to the nearest point with even coordinates. Also as in Figure \ref{f.sailing}, paths should only be allowed to fluctuate on slightly shorter time intervals which if denoted $I_1,I_2,I_3,I_4$ will satisfy $|I_1|=a+\tfrac{b-a} 2$, $|I_2|=\tfrac{b-a} 2-1$, $|I_3|=\tfrac{c-b}2 -1$ and $|I_4|=\tfrac{c-b}2 -1$. By concatenating the measure $\mu_{|I_1|},\ldots,\mu_{|I_4|}$ on directed paths in each 4 directed cones, we can defined probability measures $\mu_{x,y}$ on paths from $x$ to $y$ which satisfy the EIT property from Theorem \ref{th.UP} uniformly in $x$ and $y$ and with slightly worse constants $C,\alpha$ than in Theorem \ref{th.UP}. 
The (non-optimal) geometric condition $B(\tfrac{x+y} 2, 2\|x-y\|_2)\subset \Lambda$ is there only to ensure that the above cones remain inside $\Lambda$. \qed 

\medskip

%
%
%
\begin{remark}\label{r.Dproof}
As we pointed out earlier, the spontaneous magnetization with Dirichlet conditions (estimate~\eqref{e.onepoint}) is proved in the same way by setting $h_j \to \infty$ on the boundary and by using a cone of directed paths which connect $0$ to a point of the boundary $\p \Lambda_n$ as in Figure \ref{f.sailing}.  
\end{remark}
\smallskip

To conclude the Proof of Theorem \ref{th.disXY}, we are left with the proof of Lemma \ref{l.R}.

\smallskip
\ni
{\em Proof of Lemma \ref{l.R}.}
(N.B. This is the same proof as in \cite{abbe2018group} except we make it quantitative in $\beta$). 

From the definition of $\mathbf{R}(\omega)$ in~\eqref{e.R}, we have 
\begin{align*}\label{}
\EFK{\beta}{XY}{\mathbf{R}(\omega) \mathbf{\bar R}(\omega)}
&= \mu_n\times \mu_n \Big[ \left(\frac 1 {\lambda^2} \right)^{|p_1\cap p_2|}\Big] \\
& \hskip -3 cm \leq \sum_{k} \mu_n\times \mu_n(|p_1\cap p_2|=k) (1+ \tfrac 2 \beta)^k  \,\text{  (using~\eqref{e.lambda2} and $\beta$ sufficiently large)}\\
&\hskip -3 cm \leq  \mu_n\times \mu_n(|p_1\cap p_2| \leq  a_1 \log \beta) (1+\tfrac 2 \beta)^{a_1 \log \beta}+ C \sum_{k\geq a_1 \log \beta} e^{-\alpha k } (1 + \frac 2 \beta)^k  \\
& \hskip -3 cm \leq (1+\tfrac 2 \beta)^{a_1 \log \beta} + C \sum_{k\geq a_1 \log \beta} e^{-\frac \alpha 2 k} \\
& \hskip -3 cm \leq 1 + \tfrac{a \log \beta} \beta\,,
\end{align*}
where we used several times that the constants $\beta^*(\leq \beta)$ and   $a_1(< a)$ can be chosen large enough.  We also used Theorem \ref{th.UP} for the reflected measures $\mu_n$ for the second inequality. (N.B. One may need to take different constants $C,\alpha$ so that Theorem \ref{th.UP} applies to all measures $\mu_n$). 
\qed

\begin{remark}\label{} 
Note that the cosine function can be replaced by any other periodic functions with minor changes in the proof. In particular, the Villain model in Nishimori random disorder may be addressed the same way.
\end{remark}

\section{Symmetry breaking for Lie group valued spins}\label{s.G}

In this Section, we will adapt the above proof to the non-abelian case and prove our main Theorem \ref{th.disG}. 
See Subsection \ref{ss.G} for the relevant definitions and let $G$ be one of the compact matrix Lie groups listed in~\eqref{e.list}. 
As in the abelian case, the first main ingredient is the following Lemma.

\begin{lemma}\label{l.factorG}
For any finite domain $\Lambda \subset \Z^d$, any inverse temperature $\beta$ and any set of $k$ directed edges $(i_1j_1),\ldots (i_k j_k)$ which correspond to $k$ distinct unoriented edges, then under the Nishimori-line averaged quenched measure $\EFK{\beta}{G}{\<{\cdot}_{\omega,\beta,\Lambda}}$, the  $G$-valued random variables 
\[
\{U^*(i_m)\Omega_{i_m,j_m} U(j_m)\}_{1\leq m \leq k}
\] 
are i.i.d with  distribution $\rho_{\beta}$ on the group $G$ (see~\eqref{e.rhoG}). 
\end{lemma}

Its proof is identical to the proof of Lemma \ref{l.factor} except the {\em gauge transformations} are now given for any fixed deterministic $G$-valued field $\{S(i)\}_{i\in \Lambda}$ in $G^\Lambda$ by the following change of variable: 
\begin{align}\label{e.GT-G}
\begin{cases}
& U(i) \rightarrow S^*(i)U(i) \\
& \Omega_{i,j} \rightarrow S^*(i) \Omega_{i,j} S(j)\,. 
\end{cases}
\end{align}
Indeed, the Gibbs weight 
\begin{align*}\label{}
e^{\frac \beta 2 \sum_{(i,j)\in \vec E(\Lambda)} \Re(\Tr U^*(i) \Omega_{i,j} U(j) )} 
\end{align*}
is invariant under this transformation and by averaging uniformly over the choice of $\{S(i)\}_{i\in \Lambda}$ the partition function  cancels out exactly as in the abelian case.

\medskip

Now, as in the abelian case, one easily extracts out of Lemma \ref{l.factorG} an exact expression for the (averaged) quenched internal energy of the disordered $G$-valued spin models on the Nishimori line. To our knowledge, the internal energy of such non-abelian continuous spin systems had not been looked at before in the literature, hence we state it for the record as a Corollary (its proof is straightforward given Lemma \ref{l.factorG}). 

\begin{corollary}
Let $G$ be one of the groups listed in~\eqref{e.list}.
For any disorder $\omega$, the quenched internal energy is defined as 
\[
\<{E}_{\omega,\beta,\Lambda}:= - \frac{\p \log Z_{\omega,\beta,\Lambda}}{\p \beta}\,.
\]
On the Nishimori line $u=\beta$, the averaged quenched internal energy can be computed exactly and is given by  
\begin{align*}\label{}
\EFK{\beta}{G}{\<{E}_{\omega,\beta,\Lambda}} = -\frac {\int_G  \Re( \Tr \Omega ) e^{\beta \Re( \Tr \Omega ) } d \mu_G(\Omega)} {Z_{\rho_\beta}} |E(\Lambda)|\,. 
\end{align*}
\end{corollary}
 
\ni
Now, as in the abelian case, define for any $(i,j)\in \vec E(\Lambda)$ 
\begin{align}\label{e.12}
Y_{i,j}:= U^*(i) \Omega_{i,j} U(j)\,.
\end{align}
Lemma \ref{l.factorG} then implies that for any non-intersecting path $p$ from $x$ to $y$, 
\begin{align}\label{e.identG}
& \EFK{\beta}{G}{\<{\Re( \Tr \prod_{(ij)\in p} Y_{i,j})}_{\omega,\beta,\Lambda}} \nn  \\
& = \Re  \left( \Tr \prod_{(ij)\in p} \EFK{\beta}{G}{\<{Y_{i,j}}_{\omega,\beta,\Lambda}} \right) \nn \\
& = \Re \left( \Tr  \big(\EFK{\rho_\beta}{}{\Omega}\big)^{|p|}\right)\,.
\end{align}

Following \cite{abbe2018group}, we  define $\lambda:=\lambda_G(\beta)$ so that 
\begin{align}\label{e.lG}
\EFK{\rho_\beta}{}{\Omega} = 
\frac {\int_G \Omega \, e^{\beta \Re( \Tr \Omega ) } d \mu_G(\Omega)} {\int_G  e^{\beta \Re( \Tr \Omega ) } d \mu_G(\Omega)} = \lambda \mathbf{I}_m\,,
\end{align}
where $m$ is the trace of $\mathrm{Id}_G$ and $\mathrm{I}_m$ is the identity matrix on $\R^m$ (or $\C^m$). The fact this expectation is a multiple of $\mathrm{I}_m$ follows from the  invariance of the law $\rho_\beta$ under $G$-conjugation. In particular, this expectation needs to be in the center of $G$. Furthermore in the case where $G=U(m)$, say, one observes that $\rho_\beta$ is also invariant under complex conjugation. 
One can then check that for any group $G$ listed in~\eqref{e.list}, there exists a constant $c=c_G$ s.t. as $\beta\to \infty$, 
\[
\lambda=\lambda_G(\beta)= 1 -\tfrac{c}{\beta}+o(\beta^{-1}). 
\]
We may then rewrite~\eqref{e.identG} into the following useful identity for any $x,y\in \Lambda$ and any path $p$ from $x$ to $y$
\begin{align}\label{e.KI-G}
\EFK{\beta}{G}{\<{ \Re ( \Tr U^*(x)  \left(  \frac 1 {\lambda^{|p|}} \prod_{(ij) \in p} \Omega_{i,j} \right) U(y) }_{\omega,\beta,\Lambda} )} 
= \Tr \mathbf I_m  = m \,.
\end{align}

Let us assume from now on that $x=0$ and $y=\vec n$ (arbitrary points as well as Dirichlet boundary conditions are handled exactly as in the abelian case). 
Using the same probability measure $\mu_n$ on ``increasing'' paths from $0$ to $\vec n$, we define 
the following random matrix which is measurable w.r.t. the Nishimori disorder $\omega$:
\begin{align}\label{e.R}
\mathbf R(\omega) = \mathbf R_{0,\vec n}(\omega):=  \frac 1 {\lambda^{3n}} \EFK{\mu_n}{}{\prod_{(ij) \in p}  \Omega_{i,j}} \in  M_m(\C)\,. 
\end{align}
This random matrix satisfies by construction $\EFK{\beta}{G}{\mathbf R(\omega)} = \mathbf{I}_m$. 
Furthermore when $\beta$ is large $\mathbf R(\omega)$ is well concentrated around $\mathbf{I}_m$ in the following sense. 
\begin{lemma}[\cite{abbe2018group}]\label{l.RG}
There exist $\beta^*,a>0$ such that for all $\beta \geq \beta^*$. 
\begin{align*}\label{}
\EFK{\beta}{XY}{ \Tr \mathbf R(\omega)  \mathbf R^*(\omega)} \leq  \left(1 + \frac {a \log \beta} \beta\right) \Tr \mathbf{I}_m\,.
\end{align*}
Equivalently, 
\begin{align*}\label{}
\EFK{\beta}{XY}{ \|\mathbf R(\omega) - \mathbf{I}_m\|_F^2} \leq  \frac {a \log \beta} \beta \Tr \mathbf{I}_m\,,
\end{align*}
where $\| A\|_F= \sqrt{\Tr A^* A}$ is the Frobenius norm. 
\end{lemma}

\ni
{\em Proof.}
Following \cite{abbe2018group}, notice that 
\begin{align*}\label{}
\EFK{\beta}{XY}{\mathbf R(\omega)  \mathbf R^*(\omega)}
& =  \frac 1 {\lambda^{6n}}  \mu_n\times \mu_n \Big[ \EFK{\beta}{G}{ \prod_{(i_1j_1)\in p_1} \Omega_{i_1,j_1} \prod_{(i_2,j_2)\in \overleftarrow{p_2}} \Omega_{i_2,j_2}^*}\Big] 
\end{align*}
where the directed edges in $\overleftarrow{p_2}$ appear in the reverse order w.r.t the path $p_2\sim \mu_n$ but with same orientations as in $p_2$. Then by definition of $\lambda$ in~\eqref{e.lG}, one has  
\begin{align*}\label{}
\EFK{\beta}{XY}{\mathbf R(\omega)  \mathbf R^*(\omega)} = 
\mu_n\times \mu_n \Big[ \left(\frac 1 {\lambda^2} \right)^{|p_1\cap p_2|} \Big]
\mathbf{I}_m\,,
\end{align*} 
from which the conclusion follows exactly as in the abelian case. 
\qed

\medskip
\ni
{\em Proof of Theorem \ref{th.disG}.}
As above, we will stick to the case $x=0$ and $y=\vec n$.  Recall from~\eqref{e.KI-G} that we have
\begin{align*}\label{}
m& = \EFK{\beta}{G}{\<{\Re( \Tr U^*(x) \mathbf{R}(\omega) U(y))}_{\omega,\beta,\Lambda}}  \\
& = \EFK{\beta}{G}{\<{\Re( \Tr U^*(x)  U(y))}_{\omega,\beta,\Lambda}}  + 
\EFK{\beta}{G}{\<{\Re( \Tr U^*(x) (\mathbf{R}(\omega) - \mathbf{I}_m) U(y))}_{\omega,\beta,\Lambda}}\,.
\end{align*}
This implies
\begin{align*}\label{}
&m -  \EFK{\beta}{G}{\<{\Re( \Tr U^*(x)  U(y))}_{\omega,\beta,\Lambda}} \\
& = \EFK{\beta}{G}{\<{\Re( \Tr U^*(x) (\mathbf{R}(\omega) - \mathbf{I}_m) U(y))}_{\omega,\beta,\Lambda}} \\
& = \EFK{\beta}{G}{\<{\Re( \Tr U(y) U^*(x) (\mathbf{R}(\omega) - \mathbf{I}_m) )}_{\omega,\beta,\Lambda}}\\
&\leq
\EFK{\beta}{G}{\<{\sqrt{\Tr U(x)U^*(y)U(y) U^*(x)}  \sqrt{\Tr (\mathbf{R}^*(\omega) - \mathbf{I}_m)(\mathbf{R}(\omega) - \mathbf{I}_m)}}_{\omega,\beta,\Lambda}} \\
& = \sqrt{m} 
\EFK{\beta}{G}{\<{\sqrt{\Tr (\mathbf{R}^*(\omega) - \mathbf{I}_m)(\mathbf{R}(\omega) - \mathbf{I}_m)}}_{\omega,\beta,\Lambda}} \\
&=  
\sqrt{m} 
\EFK{\beta}{G}{\sqrt{\Tr (\mathbf{R}^*(\omega) - \mathbf{I}_m)(\mathbf{R}(\omega) - \mathbf{I}_m)}} \\
& \leq \sqrt{m} 
\EFK{\beta}{G}{\Tr (\mathbf{R}^*(\omega) \mathbf{R}(\omega) - \mathbf{I}_m)}^{1/2}
\,\,\, \text{    (since }\EFK{\beta}{G}{\mathbf R(\omega)} = \mathbf{I}_m) \\ 
& \leq \sqrt{m} \sqrt{\frac{a \log \beta}{\beta} m } = \sqrt{\frac{a \log \beta}{\beta}} \Tr \mathbf{I}_m\,,
\end{align*}
which thus concludes the proof of Theorem \ref{th.disG}. \qed 
%
%

\section{Nishimori line(s) for the classical Heisenberg model}\label{s.Heis}

In this section, we leave the setting where vertices carry group elements $U(x)\in G$ for some compact Lie group $G$ and we analyze the case of a classical Heisenberg model for a special quenched random environment. 
Let us fix a finite domain $\Lambda \subset \Z^3$. (The case of lattices $\Z^d,\, d\geq 3$ and spin $O(n)$ models with $n\geq 3$ is handled the same way). Spin configurations will be denoted
\[
S \in (\S^2)^{\Lambda}
\]

\begin{definition}[Nishimori line for classical Heisenberg model]\label{d.Heis}
Let $e_z$ be the unit vector pointing in the $z$ direction. We define  $\rho_{\beta,z}$ to be the following probability measure on $SO(3)$: 
\[
\rho_{\beta,z}(d\Omega) = \frac{ e^{\beta \<{e_z,  \Omega e_z}}}{Z_\beta} d\mu\,,
\]
where $\mu$ is the Haar measure on $SO(3)$
and $Z_\beta:= \int e^{\beta \<{e_z,  \Omega e_z}} d\mu(\Omega)$. 

The quenched disorder, which we shall still denote by $\omega$, is given by  $\omega=\{\Omega_e\}_{e\in \vec E(\Lambda)}$, where for each $e$, $\Omega_e=\Omega_e^{-1}$ and $\Omega_e$ is independent of other edges with law $\rho_{\beta,z}$. 

Given $\omega$, we consider the  following quenched classical Heisenberg model on $\Lambda$ with free boundary conditions given by the quenched Gibbs measure
\begin{align}\label{e.GibbsHeis}
dP_{\omega}(S) \propto \prod_{e=(i,j)\in \vec E(\Lambda)} \exp(\frac \beta 2  \<{S(i), \Omega_{i,j} S(j)}) \lambda^{\otimes \Lambda}(dS)\,,
\end{align}
where $\lambda$ denotes the (normalized, say) uniform measure on the sphere $\S^2$. 
As before, we will denote by $\<{\cdot}_{\omega,\beta,\Lambda}$ the expectation w.r.t this quenched Gibbs measure $P_\omega$ and $\E_{\beta}^{\Heis}$ will denote the expectation w.r.t the disorder. 
\end{definition}

\begin{remark}\label{r.notCont}
Note that as opposed to the Lie group-valued cases, for a given disorder $\omega$, the quenched Gibbs measure~\eqref{e.GibbsHeis} is not invariant under $SO(3)$. Also, if we had introduced a two parameter family $(u,\beta)$ of disorders as we did for the XY and Lie group valued spin systems, it would not be the case here that the limit $u\to \infty$ would correspond to the classical Heisenberg model. This $u\to \infty$ corresponds to an interesting model on its own which fixes the direction $e_z$ and randomly rotates transverse directions.
\end{remark}

\begin{theorem}\label{th.heis}
There exist constants $\beta^*,a >0$ s.t. for all $\beta\geq\beta^*$ then uniformly in $\Lambda\subset \Z^3$ and $x,y\in \Lambda$ sufficiently far from $\p \Lambda$, 
\begin{align}\label{e.LRH}
\EFK{\beta}{\Heis}{\<{S(x) \cdot S(y)}_{\omega,\beta,\Lambda}} \geq 1- \sqrt{\frac{a \log \beta} \beta}\,.
\end{align}
Furthermore, 
  for $\beta$ small enough,
\begin{align}\label{e.expo}
\EFK{\beta}{\Heis}{\<{S(x) \cdot S(y)}_{\omega,\beta,\Lambda}} \leq  e^{-c(\beta) \|x-y\|_2}\,.
\end{align}
\end{theorem}


\ni
{\em Proof of Theorem \ref{th.heis}.} 

Let us start by briefly explaining why despite the anisotropy of the quenched disorder, we still have a high-temperature phase. This is due to the following invariance: under free boundary conditions and under the quenched measure $P_\omega$, as we pointed out in the above remark it is no longer true (as in Subsection \ref{ss.G}, see Remark \ref{r.GlobInv}) that the system is invariant under $\{U(i)\}_{i \in \Lambda}  \mapsto \{U(i) g\}_{i\in \Lambda}$ for any global rotation $g\in G$.  Yet,  we still have the following weaker invariance under multiplication by $(-1)$, namely
\[
\{ S(i)\}_{i\in \Lambda} \mapsto \{ - S(i)\}_{i\in \Lambda}\,.
\]
Now using this invariance together with a standard Dobrushin argument at high temperature, the exponential decay in~\eqref{e.expo} easily follows. 
\medskip

We now highlight how to prove the long-range property~\eqref{e.LRH}. For simplicity, as in Section \ref{s.G}, we will stick to the case where $x=0$ and $y=\vec n$. 

\smallskip
\ni
\textbf{Step 1.} The first step of the proof is to add additional randomness to the quenched Gibbs measure $P_{\omega}$ by considering the following quenched Gibbs measure $\hat P_\omega$ on $SO(3)^\Lambda$ instead of $(\S^2)^\Lambda$:
\begin{align}\label{}
d \hat P_{\omega}(O) \propto \prod_{e=(i,j)\in \vec E(\Lambda)} \exp(\frac \beta 2  \<{O_i e_z, \Omega_{i,j} O_j e_z}) \mu^{\otimes \Lambda}(dO)\,,
\end{align}
where  $\mu$ denotes the Haar measure on $SO(3)$. It is immediate to see that if $O=\{O_i\}_{i\in \Lambda}\sim \hat P_\beta$, then $S:=\{ O_i e_z\}_{i\in \Lambda} \in (\S^2)^\Lambda$ is indeed sampled from the law $P_\omega$.  This easy step will allow us to reduce the analysis to product of group elements as in the previous sections.

Thanks to this observation, we are left with analyzing the two-point correlation
\begin{align*}\label{}
\EFK{\beta}{\Heis}{\hat E_\omega \big[ O_x e_z  \cdot O_y e_z \big]}
\end{align*}

\ni
\textbf{Step 2.} As in  previous sections, introduce the operator 
\[
\mathbf R(\omega):= \lambda^{-3n} \mathbf{E}_{\mu_n}[\prod_{e\in\gamma} \Omega_e]\,,
\]
where $\lambda=\lambda^{Heis}$ is now defined from the identity
\begin{align*}\label{}
\EFK{\rho_{\beta,z}}{}{\Omega} = \lambda \mathbf{I}_3\,. 
\end{align*}
{\em N.B. This is not the same $\lambda$ as the one  used in Section \ref{s.G} (defined in~\eqref{e.lG}) since the distribution $\rho_{\beta,z}$ is different from the distribution $\rho_\beta$ in Section \ref{s.G}.}

We then  consider the quantity 
\begin{align*}\label{}
& \EFK{\beta}{\Heis}{\hat E_\omega \big[ \<{O_x e_z , \mathbf R(\omega) O_y e_z}\big]} \\
&= \lambda^{-3n} \EFK{\beta}{\Heis}{\hat E_\omega \big[ \<{O_x e_z ,  \mathbf{E}_\mu[\prod_{e\in\gamma} \Omega_e ] \, O_y e_z}\big]} \\
& =  
\lambda^{-3n}  \mathbf{E}_\mu \Big[ 
\EFK{\beta}{\Heis}{\hat E_\omega \big[ \<{ e_z ,  O_0^{-1} \Omega_{x u_1}O_{u_1}O_{u_1}^{-1} \ldots  \Omega_{u_{3n-1}\vec n} O_{\vec n}   e_z}\big]}
\Big] 
\end{align*}

\ni
\textbf{Step 3.} By applying the same gauge transformation than in the previous sections, i.e in the present setting, for any fixed sequence of $SO(3)$ matrices $\{R_i\}_{i\in \Lambda}$, 
\begin{align*}\label{}
\begin{cases}
&O_i \to R_i^{-1} O_i \\
&\Omega_{i,j} \to R_i^{-1}  \Omega_{i,j} R_j
\end{cases}
\end{align*}
we obtain exactly as in Lemma \ref{l.factorG} that the sequence of random variables $\{O_{i}^{-1} \Omega_{i,j} O_j\}$ are i.i.d variables with law $\rho_{\beta,z}$.   This implies the identity 
\begin{align*}\label{}
\EFK{\beta}{\Heis}{\hat E_\omega \big[ \<{O_0 e_z , \mathbf R(\omega) O_{\vec n}\, e_z}\big]}
& =\<{e_z,e_z} = 1\,.
\end{align*}
As previously, we thus have 
\begin{align*}\label{}
1 - \EFK{\beta}{\Heis}{\<{S(0) \cdot S(\vec n)}_{\omega,\beta,\Lambda}} &  = 
\EFK{\beta}{\Heis}{\hat E_\omega \big[ \<{O_x e_z , (\mathbf R(\omega) - \mathbf{I}_3) O_y e_z}\big]} \\
& \leq \EFK{\beta}{\Heis}{ \| \mathbf{R}(\omega) -1\|_{op} } \\
& \leq   \EFK{\beta}{\Heis}{ \| \mathbf{R}(\omega) -1\|_{F} } \\
& \leq  \EFK{\beta}{\Heis}{ \| \mathbf{R}(\omega) -1\|_{F}^2 }^{1/2}\,.
\end{align*}
Since Lemma \ref{l.RG} also holds in the present setting (by arguing with the distribution $\rho_{\beta,z}$ instead of $\rho_\beta$), this concludes the proof of Theorem \ref{th.heis}.
\qed

\begin{remark}
In Definition \ref{d.Heis}, we may also define a Nishimori disorder by averaging over the entire $O(3)$ group rather than $SO(3)$. This gives a different model of disorder for which the proof also holds. 
\end{remark}
%


\section{Symmetry breaking of right isoclinic rotations for a Nishimori line of spin $O(4)$ model}\label{s.isoc}

As pointed out in Remark \ref{r.notCont}, the quenched Gibbs measure we introduced for the classical Heisenberg model (and by extension to all classical spin $O(n)$ models) is no longer invariant under $SO(3)$ (or under $SO(n)$ for the case of the spin $O(n)$ model). Because of this, it is only a discrete symmetry which is broken in those cases. 

The purpose of this section is to introduce a different Nishimori disorder in the special case of $\S^3$ by exploiting the fact that 
\[
\S^3 \cong  SU(2)\;\;\;\; \text{   and  } \;\;\;\;\; SO(4) \cong (SU(2) \times SU(2))/{\pm \mathrm{Id}}
\]
This allows us to rely as in Section \ref{s.G} on the above useful underlying group structure. The disorder we shall introduce below will have an interesting invariance: it will not be invariant under the whole symmetry group $SO(4)$ but it will be invariant under \textbf{right-isoclinic rotations} which correspond to $ \mathrm{Id} \times SU(2) \subset SO(4)$. 

We introduce the following identification from $\S^3$ to $SU(2)$
\[
\phi : (a,b,c,d) \in \S^3 \mapsto 
\begin{pmatrix}
a+i\cdot b & c- i\cdot d \\
- c -i \cdot d & a - i\cdot b
\end{pmatrix} \in SU(2)\,.
\]
We note that even though the map $\phi$ is not canonical (it depends on the choice of basis for $\S^4$), it satisfies the following useful identity: for any $v,w \in \S^3$, 
\[
\<{v,w} = \Re(\Tr[\phi(v)^* \phi(w)])\,.
\]
We now define a left/right actions of $SU(2)$ on $\S^3$.
\begin{definition}[Isoclinic rotations]\label{}
Any given $U\in SU(2)$ acts naturally on $\S^3$ via the following left and right actions:  for any $v\in \S^3$, 
\begin{align*}\label{}
  U\cdot v  &:= \phi^{-1}[U \phi(v)] &  \text{(left-isoclinic rotation)} \\
  v \cdot U &:= \phi^{-1}[\phi(v) U ] & \text{(right-isoclinic rotation)} 
\end{align*}
The left (resp. right) action is equivalent to the action of a {\em left-isoclinic rotation} (resp. right-isoclinic) in $SO(4)$ on $\S^3$. (N.B. the left action $SU(2)\times \mathrm{Id}   \curvearrowright \S^3$ is transitive, same for the right action). 
\end{definition}

\begin{definition}[Left-isoclinic Nishimori line for the spin $O(4)$ model]\label{d.isoc}
Let $\omega=(\Omega_{i,j})_{(ij)\in \vec E(\Lambda)} \in SU(2)^{\vec E(\Lambda)}$ be the same quenched disorder as in Section \ref{s.G} for the Lie group $G=SU(2)$. I.e. on each oriented edge $\Omega_{i,j}$ is sampled independently (subject to $\Omega_{i,j}=\Omega_{j,i}^*$) according to 
\[
d \rho_u(\Omega_{i,j}):= Z_{\rho_u}^{-1} e^{u \Re(\Tr \Omega_{i,j})} d \mu_{SU(2)}(\Omega_{i,j})\,,
\]
where $\mu_{SU(2)}$ is the Haar measure on $SU(2)$.

Given $\omega$ and using the correspondance $\phi : \S^3 \to SU(2)$, we thus define the following quenched Gibbs measure on spin configurations $\{S_i\}_{i\in \Lambda} \in (\S^3)^\Lambda$, 
\begin{align*}\label{}
\exp(\frac \beta 2  \sum_{(ij)\in \vec E} \<{S_i, \Omega_{i,j} \cdot S_j} )
&=
\exp(\frac \beta 2  \sum_{(ij)\in \vec E} \Re(\Tr[ \phi(S_i)^* \Omega_{i,j} \phi(S_j)])\,.
\end{align*}
\end{definition}
It is straightforward to check that this quenched Gibbs measure is left invariant under right-isoclinic rotations. Furthemore as $u\to \infty$, the quenched Gibbs measure converges to the classical $O(4)$ model. (This was not the case for the quenched disorder used in Section \ref{s.Heis}, see Remark \ref{r.notCont}). 

The above definition shows that this disordered model on spin $O(4)$ model can be seen as a {\em $\phi$-pull-back} of the disordered model on $G=SU(2)$ analyzed in Section \ref{s.G}. In particular, we obtain the following corollary of Theorem \ref{th.disG}.
\begin{corollary}
Let $d\geq 3$ and if $\beta$ is large enough, we have long-range order as well as  symmetry breaking of right-isoclinic rotations along the Nishimori line (i.e. $u\equiv \beta$) of the disordered spin $O(4)$ model introduced in Definition \ref{d.isoc}.
\end{corollary}


\section{Concluding remarks}\label{s.final}

\begin{remark}\label{}
Let us stress that the link between statistical reconstruction and spin systems has been very productive when spins are carried  by the vertices on trees. See for example \cite{evans2000broadcasting}. 
\end{remark}

\begin{remark}[{\bf Flexibility of this approach}]\label{r.flex}
Here are a few examples of cases which can be readily treated with the present methods.
\bnum
\item Thanks to the work by Häggström and Mossel \cite{haggstrom1998nearest}, the present long-range  order results extend to the  $2+\eps$ dimensional domains from \cite{haggstrom1998nearest}. 

\item The current setup also extends to subdomains $\Lambda \subset \Z^d$ whose edges are equipped with inhomogeneous coupling constants $(J_{i,j})_{i\sim j}$ as far as an elliptic condition $J_{i,j}\geq J >0$ is satisfied. 

Note that for classical Spin $O(n)$ model, such long-range orders with inhomogeneous elliptic coupling constants are not known due to the lack of Ginibre inequality as soon as $n\geq 3$ (\cite{ginibre1970general}).

In the same spirit, if one detects Long-Range-Order using a suitable measure on unpredictable paths, then Long-Range-Order will still hold for any graph $\hat G \supset G$ as one can use the same unpredictable path measure for $\hat G$ than for $G$. This remark allows us in particular to obtain continuous symmetry breaking for the Nishimori-line of any (reasonable) \textbf{long-range models} on $\Z^d$.  

\item  We already highlighted the fact that this approach works well with domains with arbitrary boundary $\p \Lambda$. In fact one could push this analysis further by even allowing holes in the graph as soon as one can still rely on  measures on paths $\mu_n$ which satisfy the EIT property (recall Theorem \ref{th.UP}).

\enum
\end{remark}


\begin{remark}\label{}
In the works \cite{dyson1978phase,frohlich1978phase}, reflection positivity has been used in order to prove long-range order for several quantum spin systems such as the $3d$ Heisenberg antiferromagnet.  Yet the case of the ferromagnet $3d$ quantum Heisenberg model still remains open.  It would be of great interest to try proving long-range order for this model without relying on reflection positivity (in fact it may even be necessary, see \cite{speer1985failure}). With this in mind, it would then be very interesting to generalize the present Bayesian reconstruction techniques to the setting of quantum spin systems in a suitably chosen quenched disorder.  See \cite{morita2006gauge} which initiated such an analysis.
\end{remark}

\begin{remark}\label{}
Roland Bauerschmidt has suggested that there may be a relation between the work of Kennedy-King \cite{kennedy1986spontaneous} on the abelian Higgs gauge  theory and the approach presented here. In both cases the $XY$ model is coupled with a random field.
\end{remark}

\begin{question}
Extend our results of long-range order in a neighbourhood of the Nishimori line. In particular, can one handle the case $u\to \infty$ as $\beta\to \infty$ in order to recover long-range order for classical spin $O(n)$ model?
\end{question}

\begin{question}
Can one use such techniques to establish a BKT phase transition along the Nishimori line for the $XY$ model in $d=2$ ?
\end{question}

\begin{question}
Obtain large deviation estimates for the behaviour of the $XY$ model at low temperatures in $d\geq 3$. 
\end{question}


\begin{remark}\label{}
In the work in progress \cite{GaugeNishimori}, we will extend these techniques to the case of $U(1)$ lattice gauge theory on $\Z^4$ with quenched disorder given by the Nishimori line. Namely, we shall prove a confining/deconfining phase transition for this $U(1)$ gauge theory in quenched disorder. 
\end{remark}

\section{Appendix: Correlation inequalities of Messager et al.} 

The Messager, Miracle-Sole, Pfister inequality \cite{Messager1978correlation} for the $XY$ model extends that of Ginibre \cite{ginibre1970general}. As in section \ref{s.XY}  let
$$\langle \cos\theta_0 \cos \theta_x \rangle_{\Lambda}(\omega)$$ be the expectation in the $XY$ model with disorder $\omega$ given by $e^{\beta\sum_e \cos(\theta_e - \omega_e)}$. Here $\theta_e = \theta_j -\theta_{j'}$ for an edge $e=(j,j')$.

\medskip\noindent 
\begin{theorem} 
For all $\omega$,
\[
\langle \cos\theta_0 \cos\theta_x \rangle_{\Lambda}(0) \ge \langle \cos\theta_0 \cos \theta_x \rangle_{\Lambda}(\omega).
\]
\end{theorem}
\smallskip\noindent For completeness we present  a proof of this theorem following Messager et al. 


\medskip\noindent{\it Proof.} Define the product measure by 
$$  Z(0)^{-1}Z(\omega)^{-1}e^{\beta\sum_e [\cos(\theta_e) +\cos(\theta'_e -\omega_e)] } \prod_{j\in \Lambda}d\theta_j d\theta'_j\;. $$ 
  $Z(0), Z(\omega)$ are the partition functions of the two factors.  Note that $$\cos\theta_e +\cos(\theta'_e -\omega_e) = 2\cos (\phi_e - \omega_e/2)\cos (\phi_e' - \omega_e/2)$$ 
where $$\phi_e = (\theta'_e + \theta_e)/2\,\,, \quad \phi'_e = (\theta'_e - \theta_e)/2.$$

To prove the theorem it suffices to prove that in the product measure 
$$\langle \cos\theta_0 \cos \theta_x-\cos\theta'_0 \cos \theta'_x\rangle \ge 0 \,.$$
Since $\langle \cos(\theta_0)\rangle =\langle \cos(\theta_x)\rangle=0$ the above expectation equals
$$\langle (\cos \theta_0-\cos\theta'_0)(\cos \theta_x+\cos \theta'_x)\rangle =
4\langle \sin \phi_0\cos\phi_x\sin \phi'_0\cos\phi'_x\rangle .$$
By expanding the interaction exponent in a power series in the $\phi, \phi'$ variables we see that each of the resulting terms is a square composed of identical $\phi$ and $\phi'$ factors, thus the theorem follows.
\qed

\bibliographystyle{alpha} 
\bibliography{bib-nishimori}

\newcommand{\etalchar}[1]{$^{#1}$}
\begin{thebibliography}{GHLDB85}

\bibitem[AMM{\etalchar{+}}18]{abbe2018group}
Emmanuel Abbe, Laurent Massoulie, Andrea Montanari, Allan Sly, and Nikhil
  Srivastava.
\newblock Group synchronization on grids.
\newblock {\em Mathematical Statistics and Learning}, 1(3):227--256, 2018.

\bibitem[Bal95]{balaban1995low}
Tadeusz Balaban.
\newblock A low temperature expansion for classical ${N}$-vector models. {I}.
  {A} renormalization group flow.
\newblock {\em Communications in Mathematical Physics}, 167(1):103--154, 1995.

\bibitem[Bal96]{balaban1996low}
Tadeusz Balaban.
\newblock A low temperature expansion for classical ${N}$-vector models. {II}.
  {R}enormalization group equations.
\newblock {\em Communications in mathematical physics}, 182(3):675--721, 1996.

\bibitem[Bal98a]{balaban1998large}
Tadeusz Balaban.
\newblock The large field renormalization operation for classical ${N}$-vector
  models.
\newblock {\em Communications in mathematical physics}, 198(3):493--534, 1998.

\bibitem[Bal98b]{balaban1998low}
Tadeusz Balaban.
\newblock A low temperature expansion for classical ${N}$-vector models {III}.
  a complete inductive description, fluctuation integrals.
\newblock {\em Communications in mathematical physics}, 196(3):485--521, 1998.

\bibitem[Bis09]{biskup2009reflection}
Marek Biskup.
\newblock Reflection positivity and phase transitions in lattice spin models.
\newblock In {\em Methods of contemporary mathematical statistical physics},
  pages 1--86. Springer, 2009.

\bibitem[BPP98]{benjamini1998unpredictable}
Itai Benjamini, Robin Pemantle, and Yuval Peres.
\newblock Unpredictable paths and percolation.
\newblock {\em Annals of probability}, 26(3):1198--1211, 1998.

\bibitem[DLS78]{dyson1978phase}
Freeman~J Dyson, Elliott~H Lieb, and Barry Simon.
\newblock Phase transitions in quantum spin systems with isotropic and
  nonisotropic interactions.
\newblock In {\em Statistical Mechanics}, pages 163--211. Springer, 1978.

\bibitem[EKPS00]{evans2000broadcasting}
William Evans, Claire Kenyon, Yuval Peres, and Leonard~J Schulman.
\newblock Broadcasting on trees and the {I}sing model.
\newblock {\em Annals of Applied Probability}, pages 410--433, 2000.

\bibitem[FL78]{frohlich1978phase}
J{\"u}rg Fr{\"o}hlich and Elliott~H Lieb.
\newblock Phase transitions in anisotropic lattice spin systems.
\newblock In {\em Statistical Mechanics}, pages 127--161. Springer, 1978.

\bibitem[FS82]{frohlich1982massless}
J{\"u}rg Fr{\"o}hlich and Thomas Spencer.
\newblock Massless phases and symmetry restoration in abelian gauge theories
  and spin systems.
\newblock {\em Communications in Mathematical Physics}, 83(3):411--454, 1982.

\bibitem[FSS76]{frohlich1976infrared}
J{\"u}rg Fr{\"o}hlich, Barry Simon, and Thomas Spencer.
\newblock Infrared bounds, phase transitions and continuous symmetry breaking.
\newblock {\em Communications in Mathematical Physics}, 50(1):79--95, 1976.

\bibitem[GHLDB85]{georges1985exact}
Antoine Georges, David Hansel, Pierre Le~Doussal, and J-P Bouchaud.
\newblock Exact properties of spin glasses. {II}. {N}ishimori's line: new
  results and physical implications.
\newblock {\em Journal de Physique}, 46(11):1827--1836, 1985.

\bibitem[Gin70]{ginibre1970general}
Jean Ginibre.
\newblock General formulation of {G}riffiths' inequalities.
\newblock {\em Communications in mathematical physics}, 16(4):310--328, 1970.

\bibitem[GS20]{garban2020statistical}
Christophe Garban and Avelio Sep{\'u}lveda.
\newblock Statistical reconstruction of the {G}aussian free field and {KT}
  transition.
\newblock {\em arXiv preprint arXiv:2002.12284}, 2020.

\bibitem[GS21]{GaugeNishimori}
Christophe Garban and Thomas Spencer.
\newblock Bayesian statistics and deconfining transition for ${U}(1)$ lattice
  gauge theory on the {N}ishimori line.
\newblock In preparation, 2021.

\bibitem[Gut80]{guth1980existence}
Alan~H Guth.
\newblock Existence proof of a nonconfining phase in four-dimensional u (1)
  lattice gauge theory.
\newblock {\em Physical Review D}, 21(8):2291, 1980.

\bibitem[HM98]{haggstrom1998nearest}
Olle H{\"a}ggstr{\"o}m and Elchanan Mossel.
\newblock Nearest-neighbor walks with low predictability profile and
  percolation in $2+\epsilon$ dimensions.
\newblock {\em Annals of probability}, 26(3):1212--1231, 1998.

\bibitem[Hof98]{hoffman1998unpredictable}
Christopher Hoffman.
\newblock Unpredictable nearest neighbor processes.
\newblock {\em The Annals of Probability}, 26(4):1781--1787, 1998.

\bibitem[Iba99]{iba1999nishimori}
Yukito Iba.
\newblock The {N}ishimori line and {B}ayesian statistics.
\newblock {\em Journal of Physics A: Mathematical and General}, 32(21):3875,
  1999.

\bibitem[JP02]{jacobsen2002phase}
Jesper~Lykke Jacobsen and Marco Picco.
\newblock Phase diagram and critical exponents of a potts gauge glass.
\newblock {\em Physical Review E}, 65(2):026113, 2002.

\bibitem[KK86]{kennedy1986spontaneous}
Tom Kennedy and Chris King.
\newblock Spontaneous symmetry breakdown in the abelian higgs model.
\newblock {\em Communications in mathematical physics}, 104(2):327--347, 1986.

\bibitem[KT73]{kosterlitz1973ordering}
John~Michael Kosterlitz and David~James Thouless.
\newblock Ordering, metastability and phase transitions in two-dimensional
  systems.
\newblock {\em Journal of Physics C: Solid State Physics}, 6(7):1181, 1973.

\bibitem[LDH88]{le1988location}
Pierre Le~Doussal and A~Brooks Harris.
\newblock Location of the ising spin-glass multicritical point on nishimori's
  line.
\newblock {\em Physical review letters}, 61(5):625, 1988.

\bibitem[MMSP78]{Messager1978correlation}
A~Messager, S~Miracle-Sole, and Ch~Pfister.
\newblock Correlation inequalities and uniqueness of the equilibrium state for
  the plane rotator ferromagnetic model.
\newblock {\em Communications in Mathematical Physics}, 58(1):19--29, 1978.

\bibitem[MON06]{morita2006gauge}
Satoshi Morita, Yukiyasu Ozeki, and Hidetoshi Nishimori.
\newblock Gauge theory for quantum spin glasses.
\newblock {\em Journal of the Physical Society of Japan}, 75(1):014001--014001,
  2006.

\bibitem[MW66]{mermin1966absence}
David Mermin and Herbert Wagner.
\newblock Absence of ferromagnetism or antiferromagnetism in one-or
  two-dimensional isotropic heisenberg models.
\newblock {\em Physical Review Letters}, 17(22):1133, 1966.

\bibitem[Nis81]{nishimori1981internal}
Hidetoshi Nishimori.
\newblock Internal energy, specific heat and correlation function of the
  bond-random ising model.
\newblock {\em Progress of Theoretical Physics}, 66(4):1169--1181, 1981.

\bibitem[Nis01]{nishimori2001statistical}
Hidetoshi Nishimori.
\newblock {\em Statistical physics of spin glasses and information processing:
  an introduction}.
\newblock Number 111. Clarendon Press, 2001.

\bibitem[Nis02]{nishimori2002exact}
Hidetoshi Nishimori.
\newblock Exact results on spin glass models.
\newblock {\em Physica A: Statistical Mechanics and its Applications},
  306:68--75, 2002.

\bibitem[ON93]{ozeki1993phase}
Yukiyasu Ozeki and Hidetoshi Nishimori.
\newblock Phase diagram of gauge glasses.
\newblock {\em Journal of Physics A: Mathematical and General}, 26(14):3399,
  1993.

\bibitem[Sin11]{singer2011angular}
Amit Singer.
\newblock Angular synchronization by eigenvectors and semidefinite programming.
\newblock {\em Applied and computational harmonic analysis}, 30(1):20--36,
  2011.

\bibitem[Spe85]{speer1985failure}
Eugene~R Speer.
\newblock Failure of reflection positivity in the quantum heisenberg
  ferromagnet.
\newblock {\em letters in mathematical physics}, 10(1):41--47, 1985.

\bibitem[Tan02]{tanaka2002statistical}
Kazuyuki Tanaka.
\newblock Statistical-mechanical approach to image processing.
\newblock {\em Journal of Physics A: Mathematical and General}, 35(37):R81,
  2002.

\bibitem[WS13]{wang2013exact}
Lanhui Wang and Amit Singer.
\newblock Exact and stable recovery of rotations for robust synchronization.
\newblock {\em Information and Inference: A Journal of the IMA}, 2(2):145--193,
  2013.

\end{thebibliography}
\end{document}